\documentclass[11pt,a4paper]{amsart}

\usepackage[margin=1.15in]{geometry}
\usepackage[T1]{fontenc}
\usepackage{lmodern}
\usepackage{amsmath,amssymb,amsthm}
\usepackage{enumitem,needspace}
\usepackage{hyperref}
\hypersetup{pdftitle={On Chern conjecture for complete minimal hypersurfaces with constant scalar curvature in spheres}}

\numberwithin{equation}{section}
\allowdisplaybreaks[2]
\newtheorem*{cconj}{Chern Conjecture}
\newtheorem{thm}{Theorem}[section]
\newtheorem{cor}[thm]{Corollary}
\newtheorem{lem}[thm]{Lemma}
\theoremstyle{remark}
\newtheorem*{rem}{Remark}

\newcommand{\Sph}{\mathbb S}
\newcommand{\tr}{\operatorname{tr}}

\title[On Chern conjecture for complete minimal hypersurfaces]{On Chern conjecture for complete minimal hypersurfaces with constant scalar curvature in spheres}
\author{Juan Li}
\address{School of Mathematics, Hangzhou Normal University, Hangzhou, China}
\email{juanli@zju.edu.cn}

\author{Hongwei Xu}
\address{Center of Mathematical Sciences, Zhejiang University, Hangzhou, China}
\email{xuhw@zju.edu.cn}

\author{Zhiyuan Xu}
\address{School of Mathematics, Hangzhou Normal University, Hangzhou, China}
\email{xuzy@hznu.edu.cn}

\author{Entao Zhao}
\address{Center of Mathematical Sciences, Zhejiang University, Hangzhou, China}
\email{zhaoet@zju.edu.cn}

\dedicatory{Dedicated to Professor Yuanlong Xin on the occasion of his 85th birthday}

\subjclass[2020]{Primary 53C42; Secondary 53C24, 53C40.}
\keywords{Chern conjecture, the second gap problem, constant scalar curvature, minimal hypersurfaces.}
\thanks{This research was partly supported by the National Natural Science Foundation of China, Grant Nos. 12471051,12171423.}

\begin{document}

\begin{abstract}
Let $x:M^n\to\Sph^{n+1}(1)$ be a minimal immersion of a complete manifold, where $n\ge4$.
We prove that if the squared norm $S$ of the second fundamental form is constant and $S>n$, then $S>\frac{8}{5}n$. As a consequence, we obtain the same gap for complete connected minimal immersions $M^n\to\Sph^{n+q}(1)$ with flat normal bundle, where $n\ge4$ and $q\ge2$.
\end{abstract}
\maketitle

\section{Introduction}

For a closed connected minimal immersion $x:M^n\to\Sph^{n+q}(1)$, let $S=|h|^2$ denote 
the squared norm of the second fundamental form. Simons~\cite{Simons} proved that 
if $0\le S\le \frac{n}{2-\frac{1}{q}}$ everywhere on $M$, then either $S\equiv0$ and $M$ is totally geodesic, or $S\equiv \frac{n}{2-\frac{1}{q}}$. In the case $S\equiv \frac{n}{2-\frac{1}{q}}$, Chern--do Carmo--Kobayashi~\cite{CdCK} identified the images as the minimal Clifford hypersurfaces for $q=1$ and the Veronese surface in $\Sph^4(1)$ for $n=q=2$. Lawson~\cite{Lawson} independently obtained the same classification for $q=1$. Thus, in codimension one, the only possible constant values of $S$ in $[0,n]$ are $0$ and $n$. By the Gauss equation, the scalar curvature of $M$ is $n(n-1)-S$, so the same gap can be expressed intrinsically. Chern~\cite{Chern1968,CdCK} conjectured that, in each fixed dimension, the entire set of constant scalar curvatures realized by closed minimal hypersurfaces is discrete; the refined version and the stronger version were formulated by Verstraelen--Montiel--Ros--Urbano and Xu--Xu, respectively (see \cite{Verstraelen1986,XuXu,XuXuSurvey}).

\begin{cconj}
Let $M^n$ be a closed connected minimal hypersurface in $\Sph^{n+1}(1)$.
\begin{enumerate}[label=\textup{(\Alph*)},leftmargin=*,nosep]
\item \textup{(Standard version)} For fixed $n$, the possible constant values of the scalar curvature among such hypersurfaces form a discrete subset of $\mathbb R$.
\item \textup{(Refined version)} If $M$ has constant scalar curvature, then $M$ is isoparametric.
\item \textup{(Stronger version)} Set $a_k=(k-\operatorname{sgn}(5-k))n$ for $1\le k\le5$.
  \begin{enumerate}[label=\textup{(\roman*)},leftmargin=*,nosep]
  \item For each $1\le k\le4$, if $a_k\le S\le a_{k+1}$ on $M$, then $M$ is isoparametric and either $S\equiv a_k$ or $S\equiv a_{k+1}$.
  \item If $S\ge a_5$ on $M$, then $M$ is isoparametric and $S\equiv a_5$.
  \end{enumerate}
\end{enumerate}
\end{cconj}

The values $a_1,\ldots,a_5$ are $0,n,2n,3n,5n$, respectively. They arise from the formula $S=(g-1)n$ 
for minimal isoparametric hypersurfaces, where the number $g$ of distinct principal curvatures 
belongs to $\{1,2,3,4,6\}$. The interval $[a_2,a_3]$ extends from the Clifford value $n$ to the value $2n$ attained by Cartan's examples with three distinct principal curvatures. The second gap problem asks whether a closed minimal hypersurface can have constant $S$ in the open interval $(n,2n)$. In 1983, Peng--Terng~\cite{PengTerngChapter} gave the first answer to the second gap problem, proving that $S>n+\frac{1}{12n}$ for closed minimal hypersurfaces with constant $S>n$. In dimension three, they further showed that no such hypersurface can have constant $S\in(3,6)$, while Cartan's example attains the endpoint $S=6$.

In $\Sph^4(1)$, de Almeida--Brito~\cite{AlmeidaBrito} proved that a closed hypersurface with constant mean curvature and constant nonnegative scalar curvature is isoparametric. Since nonnegativity of the scalar curvature amounts to $S\le 6$ in this case, Chang~\cite{ChangJDG} removed this restriction by establishing the refined conjecture for closed minimal hypersurfaces in $\Sph^4(1)$. In $\Sph^5(1)$, Deng--Gu--Wei~\cite{DengGuWei} proved that a closed minimal hypersurface with constant scalar curvature is isoparametric if it is Willmore, a condition equivalent in this setting to the vanishing of the cubic trace $f_3:=\tr(h^3)$. He--Xu--Zhao~\cite{HeXuZhao} relaxed this condition by assuming only that $f_3$ is constant and obtained the same conclusion. See \cite{DengKou,GLLY} for more relevant results. 

For $n\ge4$, Yang--Cheng~\cite{YangCheng} improved the pinching constant from $\frac{1}{12n}$ to $\frac{n}{3}$, proving that every closed minimal hypersurface with constant $S>n$ satisfies $S\ge\frac{4}{3}n$. For closed minimal hypersurfaces with constant $S>n$, Suh--Yang~\cite{SuhYang} proved $S>\frac{10}{7}n$. Recently, Chen~\cite{Chen} sharpened this bound to $S>\frac{85}{59}n$ in the closed case; Cheng--Li--Wei--Yin~\cite{CLWY} obtained $S>\frac{3}{2}n$ in the complete case. If $f_3$ is also constant, Cheng--Wei--Yamashiro~\cite{CWY} obtained a stronger bound, and Ge--Tan--Yan--Zhang~\cite{GTYZ} as well as Guan~\cite{Guan} obtained a sharp bound by proving $S>n$ implies $S\ge2n$. Assuming embeddedness, Tan--Tang--Xie--Yan \cite{TanTangetc} proved that the set of scalar curvature values attained by closed minimal hypersurfaces with constant  $S$ and $f_3$ is locally finite.

If we do not require that the scalar curvature is constant, then the problem becomes more difficult. Peng--Terng~\cite{PengTerngAnn} provided the first answer to the second pinching problem, proving that there exists a positive constant $\delta(n)$ depending only on $n$, such that if $n \leqslant S \leqslant n + \delta(n)$, $n \leqslant 5$, then $S \equiv n$. 
Later, Wei--Xu~\cite{WeiXu} extended the result to $n = 6, 7$ and Zhang~\cite{Q.Zhang} promoted it to $n \leqslant 8$.  Finally, Ding--Xin~\cite{DingXin} extended it to all dimensions by establishing a new estimate; in particular, they showed that if the dimension is $n \geqslant 6$, then the pinching constant $\delta(n) = \frac{n}{23}$. After that, Xu--Xu~\cite{XuXu} improved it to $\delta(n) = \frac{n}{22}$ and Lei--Xu--Xu~\cite{LeiXuXu} showed $\delta(n) = \frac{n}{18}$.

In this paper, we mainly prove the following second gap theorem for complete minimal hypersurfaces in the unit sphere.
\begin{thm}\label{main}
Let $x:M^n\to\Sph^{n+1}(1)$ be a minimal immersion, where $n\ge4$ and $M^n$ is complete in the induced metric. If the squared norm $S$ of the second fundamental form is constant and $S>n$, then $S>\frac{8}{5}n$.
\end{thm}

For a complete connected minimal immersion $M^n\to\Sph^{n+q}(1)$ with flat normal bundle, where $n\ge3$ and $q\ge2$, Li--Xu~\cite{LiXu} proved that $S<2n$ forces the image to lie in a fixed totally geodesic subsphere $\Sph^{n+1}$. Combining their result with Theorem~\ref{main} gives the following corollary.

\begin{cor}\label{cor:flat-normal}
Let $x:M^n\to\Sph^{n+q}(1)$ be a minimal immersion with flat normal bundle, where $n\ge4$, $q\ge2$ and $M^n$ is connected and complete in the induced metric. If $S$ is constant and $S>n$, then $S>\frac{8}{5}n$.
\end{cor}

To prove Theorem~\ref{main}, suppose $S$ is constant and $n<S\le\frac{8n}{5}$. In Section~\ref{sec:pre} we apply the Omori--Yau maximum principle to a polynomial $F$ in $f_3^2$ and $\tr(h^4)$. Along a maximizing sequence, $|\nabla F|\to0$ gives an orthogonality relation for the limits of the diagonal entries $(h_{iik})_i$ in principal frames. We use it to estimate $\nabla f_3$ and the commutators $[h,\nabla_k h]$. Projecting the fully symmetric part of $\nabla^2h$ onto two orthogonal tensors bounds its squared norm from below. In Section~\ref{sec:scalar-inequality}, these estimates, the Peng--Terng identity and $\limsup\Delta F\le0$ yield two scalar inequalities. We sharpen the first by retaining the exact projection denominator $\frac{S^2+2\tr(h^4)}{3}$, the squared norm of the symmetrization of $h\otimes h$. In their weighted sum we fix the componentwise limits of $\frac{h}{\sqrt S}$ and $\frac{\nabla h}{|\nabla h|}$ and rescale to $|h|^2=s$, $|\nabla h|^2=s(s-n)$ by Simons' identity. Section~\ref{sec:spectral-reduction} chooses the weights using only $n$ and $\frac{f_3^2}{S^3}$. For these fixed weights the right-hand side of the weighted inequality minus its left-hand side is nonnegative at $s=S$ and strictly increasing for $s>n$. It must therefore be nonnegative at $s=\frac{8n}{5}$, contradicting the strictly negative upper bound in Lemma~\ref{lem:endpoint-verification}.

\section{Preliminaries}\label{sec:pre}

Let $x:M^n\to\Sph^{n+1}(1)$ be a complete minimal immersion with $n\ge4$, and $g$ be the induced metric on $M$. The indices $i,j,k,\ell,\ldots$ range from $1$ to $n$, and $A,B,C,\ldots$ range from $1$ to $n+1$.
Choose a local orthonormal frame $\{e_A\}$ along $M$ with $e_1,\ldots,e_n$ tangent to $M$ and $e_{n+1}=\nu$ normal to $M$. The dual coframe $\{\theta_A\}$ satisfies $\theta_{n+1}|_{TM}=0$. The connection forms are $\theta_{AB}(X)=\langle\bar\nabla_Xe_A,e_B\rangle$, where $\bar\nabla$ is the Levi--Civita connection of the sphere.

The second fundamental form of $M$ in $\Sph^{n+1}(1)$ is $h\,\nu$, where $h(X,Y)=\langle\bar\nabla_XY,\nu\rangle$ for tangent vector fields $X,Y$. In this frame,
\[
 h=\sum_{i,j=1}^n h_{ij}\,\theta_i\otimes\theta_j,\qquad
 h_{ij}=h(e_i,e_j)=h_{ji},\qquad
 \theta_{i,n+1}=\sum_{j=1}^n h_{ij}\,\theta_j .
\]
We use the induced metric to identify $h$ with the shape operator. The mean curvature vector $\mathbf H$ and the squared norm $S$ of the second fundamental form are given by
\[
 \mathbf H=(\tr_g h)\nu=\left(\sum_{i=1}^n h_{ii}\right)\nu,\qquad
 S=|h|^2=\sum_{i,j=1}^n h_{ij}^2 .
\]
Since $M$ is minimal, $\mathbf H=0$ and $\sum_i h_{ii}=0$. Throughout the proof, $S$ is a constant greater than $n$.

Let $\nabla$ be the Levi--Civita connection of $g$ and set $\nabla_k=\nabla_{e_k}$. Together with the normal connection $\nabla^\perp$, it induces a connection on normal-valued tensors, also denoted by $\nabla$. Since $\nabla^\perp\nu=0$, we have $\nabla^r(h\,\nu)=(\nabla^r h)\,\nu$ for $r\ge1$. The components of $\nabla h$ and $\nabla^2h$ are defined by
\[
\begin{aligned}
 \sum_k h_{ijk}\,\theta_k
 &=\mathrm{d}h_{ij}+\sum_m\bigl(h_{mj}\,\theta_{mi}+h_{im}\,\theta_{mj}\bigr),\\
 \sum_\ell h_{ijk\ell}\,\theta_\ell
 &=\mathrm{d}h_{ijk}+\sum_m\bigl(h_{mjk}\,\theta_{mi}+h_{imk}\,\theta_{mj}
                                  +h_{ijm}\,\theta_{mk}\bigr).
\end{aligned}
\]
Summing over ordered indices for the induced inner product $\langle\cdot,\cdot\rangle$ gives
\[
 |\nabla h|^2=\sum_{i,j,k}h_{ijk}^2,\qquad
 |\nabla^2h|^2=\sum_{i,j,k,\ell}h_{ijk\ell}^2 .
\]
We take $\Delta=\tr_g\nabla^2$, so $(\Delta h)_{ij}=\sum_k h_{ijkk}$ and $\Delta f=\operatorname{div}(\nabla f)$ for functions. Since $S$ is constant, Simons' identity $\frac{1}{2}\Delta S=|\nabla h|^2+(n-S)S$ gives
\begin{equation}\label{eq:simons-gradient}
 |\nabla h|^2=S(S-n)>0 .
\end{equation}

At a point $p\in M$, choose a principal frame with $h_{ij}=\lambda_i\delta_{ij}$, where $\sum_i\lambda_i=0$ and $\sum_i\lambda_i^2=S$. Set $\mu_i=\frac{\lambda_i}{\sqrt{S}}$ and $f_j=\tr(h^j)=\sum_i\lambda_i^j$.
The Codazzi equation and the symmetry of the second fundamental form give
\[
 h_{ijk}=h_{ikj}=h_{jik}.
\]

The curvature tensor $R$ of $g$ has components $R_{ijk\ell}=\langle R(e_k,e_\ell)e_j,e_i\rangle$. The Ricci commutation formula and the Gauss equation for a hypersurface in the unit sphere are
\[
\begin{aligned}
 h_{ijk\ell}-h_{ij\ell k}
 &=\sum_m\bigl(h_{mj}\,R_{mik\ell}+h_{im}\,R_{mjk\ell}\bigr),\\
 R_{ijk\ell}
 &=\delta_{ik}\delta_{j\ell}-\delta_{i\ell}\delta_{jk}
   +h_{ik}h_{j\ell}-h_{i\ell}h_{jk}.
\end{aligned}
\]
The Gauss equation gives $R_{ikik}=1+\lambda_i\lambda_k$ for $i\ne k$ and
\[
 \operatorname{Ric}
 =\sum_{i,j,k}R_{ikjk}\,\theta_i\otimes\theta_j
 =\sum_i(n-1-\lambda_i^2)\,\theta_i\otimes\theta_i
 \ge(n-1-S)g .
\]
Differentiating the Codazzi equation and then commuting derivatives gives
\begin{equation}\label{eq:ricci-commutation}
 h_{iikk}-h_{kkii}=h_{ikik}-h_{ikki}=(\lambda_i-\lambda_k)R_{ikik}
 =(\lambda_i-\lambda_k)(1+\lambda_i\lambda_k).
\end{equation}

Set
\[
 \mathcal{A}=\sum_{ijk}\lambda_i^2h_{ijk}^2,\qquad
 \mathcal{B}=\sum_{ijk}\lambda_i\lambda_jh_{ijk}^2,\qquad
 \mathcal{C}=\sum_{ijk}\lambda_i h_{ijk}^2.
\]
Since $|\nabla h|^2$ is constant, the Peng--Terng identity gives
\begin{equation}\label{eq:derivative-bounds}
 |\nabla^2h|^2=(S-2n-3)|\nabla h|^2+3(\mathcal{A}-2\mathcal{B}).
\end{equation}
Summing the squares of the commutator components $[h,\nabla_k h]_{ij}=(\lambda_i-\lambda_j)h_{ijk}$ gives
\begin{equation}\label{eq:shape-commutator}
 \sum_k|[h,\nabla_k h]|^2=\sum_{ijk}(\lambda_i-\lambda_j)^2h_{ijk}^2=2(\mathcal{A}-\mathcal{B}) .
\end{equation}

Cyclicity of the trace gives
\[
 \nabla_k f_j=j\sum_i\lambda_i^{j-1}h_{iik}.
\]
Applying the product rule and using $\Delta h=(n-S)h$, we obtain
\begin{equation}\label{eq:trace-laplacians}
 \frac{\Delta f_3}{3}=(n-S)f_3+2\mathcal{C},\qquad
 \frac{\Delta f_4}{4}=(n-S)f_4+2\mathcal{A}+\mathcal{B}.
\end{equation}

Differentiating $\tr h=0$ and $|h|^2=S$ along $e_k$, we obtain
\[
 \langle\nabla_k h,I\rangle=\sum_i h_{iik}=0,\qquad
 \langle\nabla_k h,h\rangle=\sum_i\lambda_i h_{iik}=0,
\]
where $I$ is the identity on $T_pM$.

Since $I\perp h$, the orthogonal projection $K$ of $\frac{h^2}{S}$ onto $\{I,h\}^{\perp}$ satisfies
\[
 K=\frac{h^2}{S}-\frac{f_3}{S^2}h-\frac{I}{n},\qquad
 |K|^2=\frac{f_4}{S^2}-\frac{f_3^2}{S^3}-\frac{1}{n}.
\]
For $j=3$, the two trace identities give
\[
 \nabla_k f_3=3\langle h^2,\nabla_k h\rangle=3S\langle K,\nabla_k h\rangle.
\]
Thus $K=0$ exactly at points with at most two distinct principal curvatures.

To obtain a third relation for the limits of $(h_{iik})_i$ along a maximizing sequence, consider the function
\[
\begin{aligned}
 F=\frac{Sf_4}{4}-\frac{137f_3^2}{600}+\frac{f_3^4}{120S^3}=\frac{S^3}{4}\left(|K|^2+\frac{1}{n}\right)+\frac{13f_3^2}{600}+\frac{f_3^4}{120S^3}.
\end{aligned}
\]
Both $f_4$ and $f_3^2$ are unchanged on reversing $\nu$, so $F$ is globally defined. Since $S$ is constant and $f_3^2\le Sf_4\le S^3$, $F$ is bounded above. As $M$ is complete and $\operatorname{Ric}\ge(n-1-S)g$, the Omori--Yau maximum principle~\cite{Omori,Yau} gives a sequence $(p_j)$ in $M$ such that
\[
 F(p_j)\longrightarrow\sup_M F,\qquad
 |\nabla F|(p_j)\longrightarrow0,\qquad
 \limsup_{j\to\infty}\Delta F(p_j)\le0 .
\]

At each $p_j$, choose a unit normal with $f_3(p_j)\ge0$ and an orthonormal principal frame. Combining \eqref{eq:simons-gradient} and \eqref{eq:derivative-bounds} with $|\mathcal A|,|\mathcal B|\le S|\nabla h|^2$, we obtain a uniform bound for $|\nabla^2h|$. Pass to a subsequence, still denoted by $(p_j)$, along which all components of $h,\nabla h,\nabla^2h$ in these frames converge. We keep the notation $\lambda_i,h_{ijk},h_{ijk\ell}$ for their limits and use the same notation for their contractions.

Holding $S$ and $f_4$ fixed, we obtain
\[
 \partial_{f_3}F=-\frac{137 f_3}{300}+\frac{f_3^3}{30S^3},\qquad
 \partial_{f_3}^2F=-\frac{137}{300}+\frac{f_3^2}{10S^3}.
\]
The chain rule gives
\[
\begin{aligned}
 \nabla_k F&=\frac{S}{4}\nabla_k f_4+(\partial_{f_3}F)\nabla_k f_3,\\
 \Delta F&=\frac{S}{4}\Delta f_4+(\partial_{f_3}F)\Delta f_3
           +(\partial_{f_3}^2F)|\nabla f_3|^2.
\end{aligned}
\]
Using \eqref{eq:trace-laplacians} and passing to the limit along $(p_j)$ gives
\begin{equation}\label{F}
\begin{gathered}
  \sum_i\left(\lambda_i^3+\frac{3\,\partial_{f_3}F}{S}\lambda_i^2\right)h_{iik}=0,\\
  (n-S)(Sf_4+3f_3\,\partial_{f_3}F)+S(2\mathcal{A}+\mathcal{B})
    +6\mathcal{C}\,\partial_{f_3}F+(\partial_{f_3}^2F)|\nabla f_3|^2\le0,
\end{gathered}
\end{equation}
where the first identity holds for each tangent index $k$.

Identifying diagonal operators with their diagonal vectors in $\mathbb R^n$, we take the orthogonal projection with respect to the trace inner product
\[
 P\colon\mathbb R^n\longrightarrow
 \left(\operatorname{span}\left\{I,\,h,\,h^3+\frac{3\,\partial_{f_3}F}{S}h^2\right\}\right)^{\perp}.
\]
The two trace identities and \eqref{F} give $(h_{iik})_{i=1}^n\in\operatorname{im}P$ for each $k$. Since $K$ is diagonal, it follows that for every tangent vector $\eta$,
\[
 \langle\nabla f_3,\eta\rangle=3S\langle PK,\nabla_\eta h\rangle .
\]
In particular, $PK=0$ implies $\nabla f_3=0$.

The projection of the $k$th coordinate vector onto $\{I,h\}^{\perp}$ has squared norm $1-\frac{1}{n}-\mu_k^2$. Applying the Cauchy--Schwarz inequality to the inner products of $(h_{iik})_i$ with $PK$ and with this projected vector gives
\begin{equation}\label{eq:component-projection-bounds}
 (\nabla_k f_3)^2\le 9S^2|PK|^2\sum_i h_{iik}^2,\qquad h_{kkk}^2\le\left(1-\frac{1}{n}-\mu_k^2\right)\sum_i h_{iik}^2 .
\end{equation}

Cheng--Li--Wei--Yin~\cite{CLWY} proved the estimate
\begin{equation}\label{audit:lowerAB}
  \mathcal{A}+2\mathcal{B}\ge \frac{3\mathcal{C}^2}{|\nabla h|^2}+\frac{4(n+2)}{9}\sum_k\frac{(\nabla_k f_3)^2}{(n+2)S+2n\lambda_k^2}.
\end{equation}

\begin{lem}\label{mixed}
Let $x:M^n\to\Sph^{n+1}(1)$ be a minimal immersion with $n\ge4$ and constant $S>n$. At each point choose a unit normal for which $f_3\ge0$. Then
\[
 \frac{4(\mathcal{A}-\mathcal{B})}{S|\nabla h|^2}
 \le 2+\frac{2f_3 f_4}{S^2\bigl(S^{3/2}+f_3\bigr)} .
\]
Every sequence $(p_j)$ with $|\nabla F|(p_j)\to0$ satisfies
\begin{equation}\label{eq:trace-moment-L}
 \limsup_{j\to\infty}\frac{4(\mathcal A-\mathcal B)(p_j)}{S|\nabla h|^2(p_j)}\le2 .
\end{equation}
\end{lem}
\begin{proof}
By Codazzi symmetry and \eqref{eq:shape-commutator}, the mean coefficient of $h_{ijk}^2$ in $\frac{4(\mathcal A-\mathcal B)}{S}$ over the six permutations of three distinct indices is
\[
 \frac{2}{3}[(\mu_i-\mu_j)^2+(\mu_j-\mu_k)^2+(\mu_k-\mu_i)^2]
 \le2(\mu_i^2+\mu_j^2+\mu_k^2)\le2 .
\]
\Needspace{6\baselineskip}
Set $L:=2+\frac{2f_3 f_4}{S^2(S^{3/2}+f_3)}$ for the pointwise upper bound on $\frac{4(\mathcal A-\mathcal B)}{S|\nabla h|^2}$ and $L=2$ for the limsup bound along sequences $(p_j)$ with $|\nabla F|(p_j)\to0$.
Since $L\ge2$ in both cases, it suffices to prove for each $k$ that the terms with a repeated index satisfy
\begin{equation}\label{eq:repeated-index-group}
 4\sum_{i\ne k}(\mu_i-\mu_k)^2 h_{iik}^2
 \le Lh_{kkk}^2+3L\sum_{i\ne k}h_{iik}^2 .
\end{equation}
For fixed $k$, at most one coefficient $4(\mu_i-\mu_k)^2$ can exceed $3L$. Indeed, for distinct $i,j\ne k$,
\[
 \frac{4}{L}\bigl[(\mu_i-\mu_k)^2+(\mu_j-\mu_k)^2\bigr]
 \le\frac{12}{L}(\mu_i^2+\mu_j^2+\mu_k^2)\le6 .
\]
If there is no such index, \eqref{eq:repeated-index-group} holds term by term. Otherwise the exceptional index $i$ satisfies $\mu_i\mu_k<0$, and for every $j\ne i,k$ we have
\[
 0<p_L:=\frac{4(\mu_i-\mu_k)^2}{L}-3
 \le3-\frac{4(\mu_j-\mu_k)^2}{L}.
\]

Let $q_1,\ldots,q_n$ satisfy $\sum_jq_jh_{jjk}=0$ and $q_i\ne0$. By the weighted Cauchy--Schwarz inequality,
\[
\begin{aligned}
 p_Lq_i^2h_{iik}^2=p_L\left(q_kh_{kkk}+\sum_{j\ne i,k}q_jh_{jjk}\right)^2\le\left(p_Lq_k^2+\sum_{j\ne i,k}q_j^2\right)
      \left(h_{kkk}^2+p_L\sum_{j\ne i,k}h_{jjk}^2\right).
\end{aligned}
\]
Thus $p_Lq_k^2+\sum_{j\ne i,k}q_j^2\le q_i^2$ suffices for \eqref{eq:repeated-index-group}.

The trace identity $\sum_j\mu_jh_{jjk}=0$ permits the choice $q_j=\mu_j$. With $p_2=2(\mu_i-\mu_k)^2-3\ge p_L$, the normalization $\sum_j\mu_j^2=1$ gives
\[
 \mu_i^2-\sum_{j\ne i,k}\mu_j^2
 =p_2(1-\mu_k^2)+2(1-\mu_k^2+\mu_i\mu_k)^2 .
\]
When $\mu_k^2\le\frac12$, this gives $p_2\mu_k^2+\sum_{j\ne i,k}\mu_j^2\le\mu_i^2$, proving \eqref{eq:repeated-index-group} with $L=2$. Suppose now that $\mu_k^2>\frac12$. The inequalities $f_3\ge0$ and $\sum_{j\ne k}|\mu_j|^3\le(1-\mu_k^2)^{\frac{3}{2}}<|\mu_k|^3$ imply $\mu_k>0>\mu_i$. For $q_j=\mu_j$, the condition on the coefficients is equivalent to
\begin{equation}\label{Ecriterion}
 L(2\mu_i^2+4\mu_k^2-1)\ge4\mu_k^2(\mu_k-\mu_i)^2 .
\end{equation}
We have $\frac{f_3}{S^{3/2}}\ge\mu_k^3-(1-\mu_k^2)^{\frac{3}{2}}\ge2\mu_k^2-1$, where the difference in the second inequality factors as $(\mu_k-\sqrt{1-\mu_k^2})(1-\mu_k)(1-\sqrt{1-\mu_k^2})\ge0$. The pointwise choice of $L$ therefore satisfies
\[
 L-2=\frac{2f_3}{S^{3/2}+f_3}\sum_j\mu_j^4
 \ge\frac{2\mu_k^2-1}{\mu_k^2}(\mu_k^4+\mu_i^4)
 \ge\frac{2\mu_k^2-1}{\mu_k^2}\left(\mu_k^2+\mu_i^2-\frac{1}{2}\right).
\]
Since $\mu_i^2\le1-\mu_k^2<\frac12$ and $\mu_k-\mu_i\le2\mu_k$, we obtain
\[
\begin{aligned}
 &L(2\mu_i^2+4\mu_k^2-1)-4\mu_k^2(\mu_k-\mu_i)^2\\
 \ge&4\mu_k^2(\mu_k+\mu_i)^2
 -\frac{\mu_k^2+\mu_i^2-\frac{1}{2}}{\mu_k^2}(2\mu_k^2-1)(1-2\mu_i^2)\\
 \ge&(\mu_k^2-\mu_i^2)^2-(2\mu_k^2-1)(1-2\mu_i^2)
 =(1-\mu_k^2-\mu_i^2)^2\ge0 .
\end{aligned}
\]
This proves \eqref{Ecriterion} and hence the pointwise bound.

Choose a subsequence realizing the limsup in \eqref{eq:trace-moment-L}. The normalization $\sum_i\mu_i^2=1$ and \eqref{eq:simons-gradient} bound the components $\mu_i$ and $h_{ijk}$, so after passing to a further subsequence they all converge in principal frames. To establish \eqref{eq:repeated-index-group} with $L=2$ for these limits, only the case $p_2>0$ and $\mu_k^2>\frac12$ remains. Multiplying $\mu_i^2-\sum_{j\ne i,k}\mu_j^2\ge p_2(1-\mu_k^2)$ by $1-\mu_k^2$ gives
\[
 p_2(1-\mu_k^2)^2+\left(\sum_{j\ne i,k}\mu_j^2\right)^2\le\mu_i^4 .
\]
In particular, $\sqrt{\sum_{j\ne i,k}\mu_j^2}\le-\mu_i$ and $|\mu_j|\le-\mu_i$ for every $j\ne k$. It follows that
\[
 \beta=-\frac{3\,\partial_{f_3}F}{S^{3/2}}
 =\left(\frac{137}{100}-\frac{f_3^2}{10S^3}\right)\frac{f_3}{S^{3/2}}
 \ge\frac{f_3}{S^{3/2}}\ge\mu_k^3+\mu_i(1-\mu_k^2).
\]
The chain rule gives $\sum_j\mu_j^2(\mu_j-\beta)h_{jjk}=\lim S^{-\frac{5}{2}}\nabla_k F=0$. Combining this identity with $\sum_j\mu_jh_{jjk}=0$ yields $\sum_jq_jh_{jjk}=0$ for the coefficients
\[
 q_j=\mu_j^2(\mu_j-\beta)+\mu_k(\beta-\mu_k)_+\mu_j,\qquad
 q_k=\mu_k^2(\mu_k-\beta)_+\le(1-\mu_k^2)(\beta-\mu_i),
\]
where $q_i<0$ and $x_+=\max\{x,0\}$ denotes the positive part of $x$. The triangle inequality and $\sum_{j\ne i,k}\mu_j^4\le\bigl(\sum_{j\ne i,k}\mu_j^2\bigr)^2$ yield
\[
 \begin{aligned}
 \left(p_2q_k^2+\sum_{j\ne i,k}q_j^2\right)^{\frac{1}{2}}
 &\le(\beta-\mu_i)\mu_i^2
       +\mu_k(\beta-\mu_k)_+\sqrt{\sum_{j\ne i,k}\mu_j^2}\\
 &\le\mu_i^2(\beta-\mu_i)-\mu_i\mu_k(\beta-\mu_k)_+=|q_i|.
 \end{aligned}
\]
The coefficient condition now holds with $L=2$, proving \eqref{eq:trace-moment-L} and \eqref{eq:repeated-index-group}.
\end{proof}

Using the symmetry of $h_{ijk}$ and the identities $\sum_i h_{iik}=\sum_i\lambda_i h_{iik}=0$, we apply Bessel's inequality as in~\cite{CLWY} to obtain, for every tangent vector $\eta=\sum_k\eta_k e_k$,
\begin{equation}\label{audit:Tv}
 \left(1+\frac{2}{n}\right)\sum_{i,j}\left(\sum_k\eta_k h_{ijk}\right)^2
 +\frac{2}{S}\sum_{i,j}\left(\sum_k\lambda_k\eta_k h_{ijk}\right)^2
 \le |\nabla h|^2|\eta|^2.
\end{equation}

By the Codazzi equation, $h_{ijk\ell}$ is symmetric in its first three indices, so the fully symmetric part $u$ of $\nabla^2h$ has components
\[
 u_{ijk\ell}=\frac{1}{4}\bigl(h_{ijk\ell}+h_{\ell ijk}+h_{k\ell ij}+h_{jk\ell i}\bigr).
\]
Expanding the definition of $u$ and using \eqref{eq:ricci-commutation}, we obtain
\[
 |\nabla^2h|^2-|u|^2=\frac{3}{2}(Sf_4-f_3^2-2S^2+nS).
\]
Substituting \eqref{eq:derivative-bounds} into the preceding identity and using \eqref{eq:simons-gradient}, we obtain
\begin{equation}\label{Uidentity}
 |u|^2=3(\mathcal{A}-2\mathcal{B})-\frac{3}{2}S^3|K|^2+ \left(S-2n-\frac{3(S+n)}{2n}\right)|\nabla h|^2 .
\end{equation}

To bound $|u|^2$ from below, we project $u$ onto the two symmetric tensors $Y$ and $V$ given by
\[
\begin{aligned}
 Y_{ijk\ell}&=\frac{S}{3}\bigl(
 h_{ij}\delta_{k\ell}+h_{ik}\delta_{j\ell}+h_{i\ell}\delta_{jk}
 +h_{jk}\delta_{i\ell}+h_{j\ell}\delta_{ik}+h_{k\ell}\delta_{ij}\bigr),\\
 V_{ijk\ell}&=\frac{1}{3S^2}\bigl(
 \lambda_i h_{jk\ell}\nabla_i f_3+\lambda_j h_{ik\ell}\nabla_j f_3
 +\lambda_k h_{ij\ell}\nabla_k f_3+\lambda_\ell h_{ijk}\nabla_\ell f_3\bigr),
\end{aligned}
\]
where $(h(\nabla f_3))_k=\lambda_k\nabla_k f_3$ in the principal frame.

Define $v$ by
\begin{equation}\label{eq:weighted-gradient}
 v=
 \begin{cases}
 \displaystyle 1-\frac{2|h(\nabla f_3)|^2}{9S^3|PK|^2|\nabla h|^2}, & PK\ne0,\\[5pt]
 1, & PK=0.
 \end{cases}
\end{equation}
We have $0\le v\le1$.

\begin{lem}\label{lem:hessian-projection}
The fully symmetric part $u$ of $\nabla^2h$ satisfies
\begin{equation}\label{UVlower}
 \frac{|u|^2}{S^3}
 \ge\frac{3(S-n)^2}{2(n+4)S^2}
 +\frac{9(n+2)(1-v)|PK|^2}{2(n+2+3nv)} .
\end{equation}
\end{lem}
\begin{proof}
For fixed $k$, the Codazzi identity $2u_{iikk}=h_{iikk}+h_{kkii}$ and \eqref{eq:ricci-commutation} give
\[
 2\sum_i\lambda_i u_{iikk}
 =-2|\nabla_k h|^2+S\lambda_k^2-f_3\lambda_k-S.
\]
For any $E$ diagonal in the principal frame, we sum this identity against $E_{kk}$ and contract $Y$ and $V$ with $h\otimes E$ to obtain
\begin{subequations}\label{eq:weighted-traces}
\begin{align}
 \langle u,h\otimes E\rangle
 &=-\sum_k E_{kk}|\nabla_k h|^2+\frac{S^2}{2}\langle K,E\rangle+\frac{S(S-n)}{2n}\tr E,\label{eq:weighted-trace-u}\\
 \langle Y,h\otimes E\rangle
 &=\frac{S}{3}\bigl(S\tr E+4\langle h^2,E\rangle\bigr),\label{eq:weighted-trace-y}\\
 \langle V,h\otimes E\rangle
 &=\frac{2}{3S^2}\sum_i\lambda_i^2(\nabla_i f_3)\langle E,\nabla_i h\rangle .\label{eq:weighted-trace-v}
\end{align}
\end{subequations}
Since $\langle T,Y\rangle=2S\langle T,h\otimes I\rangle$ for every fully symmetric $4$-tensor $T$, taking $E=I$ gives
\[
 \langle u,Y\rangle=-S|\nabla h|^2,\qquad
 |Y|^2=\frac{2(n+4)S^3}{3},\qquad Y\perp V .
\]

If $h(\nabla f_3)=0$, then $V=0$ and $v=1$, so projection onto $Y$ gives the required bound $|u|^2\ge\frac{\langle u,Y\rangle^2}{|Y|^2}=\frac{3S(S-n)^2}{2(n+4)}$.

Assume $h(\nabla f_3)\ne0$. Since $|\nabla h|^2$ is constant, $\sum_{ijk}h_{ijk}h_{ijk\ell}=0$. Codazzi symmetry and Ricci commutation give
\[
\begin{aligned}
 4\sum_{ijk}h_{ijk}u_{ijk\ell}=3\sum_{ijk}h_{ijk}(h_{ij\ell k}-h_{ijk\ell})=-6\sum_i(\lambda_i-\lambda_\ell)(1+\lambda_i\lambda_\ell)h_{ii\ell}.
\end{aligned}
\]
The trace identities reduce this to
\[
 \sum_{i,j,k}h_{ijk}u_{ijk\ell}
 =-\frac{3}{2}\lambda_\ell\sum_i\lambda_i^2h_{ii\ell}
 =-\frac{\lambda_\ell\nabla_\ell f_3}{2}.
\]
Contracting $V$ with the symmetric tensor $u$ now gives
\[
 \langle u,V\rangle
 =-\frac{2}{3S^2}\sum_\ell(\lambda_\ell\nabla_\ell f_3)^2
 =-\frac{2|h(\nabla f_3)|^2}{3S^2}.
\]

Taking $\eta=h^2(\nabla f_3)$ in the projection formula gives
\[
 |h(\nabla f_3)|^2
 =3S\left\langle PK,\sum_k\lambda_k^2(\nabla_k f_3)\nabla_k h\right\rangle.
\]
Thus $PK\ne0$, and the Cauchy--Schwarz inequality gives
\[
 \frac{2}{S}\sum_{i,j}\left(\sum_k\lambda_k^2 h_{ijk}\nabla_k f_3\right)^2
 \ge \frac{2|h(\nabla f_3)|^4}{9S^3|PK|^2}
 =(1-v)|h(\nabla f_3)|^2|\nabla h|^2 .
\]
Taking $\eta=h(\nabla f_3)$ in \eqref{audit:Tv} and applying this lower bound to the second term gives $v\ge0$, since
\[
 \left(1+\frac{2}{n}\right)\sum_{i,j}\left(\sum_k\lambda_k h_{ijk}\nabla_k f_3\right)^2
 \le v|h(\nabla f_3)|^2|\nabla h|^2.
\]
Expanding the squared norm of $3S^2V$ gives
\begin{equation}\label{eq:symmetric-test-norm}
\begin{aligned}
 9S^4|V|^2
 &=4|h(\nabla f_3)|^2|\nabla h|^2+12\sum_{i,j}\left(\sum_k\lambda_k h_{ijk}\nabla_k f_3\right)^2\\
 &\le\frac{4(n+2+3nv)}{n+2}|h(\nabla f_3)|^2|\nabla h|^2 .
\end{aligned}
\end{equation}
Since $Y\perp V$, Bessel's inequality gives
\[
 |u|^2\ge\frac{\langle u,Y\rangle^2}{|Y|^2}+\frac{\langle u,V\rangle^2}{|V|^2}
 \ge\frac{3S(S-n)^2}{2(n+4)}+\frac{(n+2)|h(\nabla f_3)|^2}{(n+2+3nv)|\nabla h|^2}.
\]
Dividing by $S^3$ and using \eqref{eq:weighted-gradient} gives \eqref{UVlower}.
\end{proof}

For the contractions with $h\otimes h$ and $h\otimes K$, set
\[
\begin{aligned}
 Z_{ijk\ell}
 &=\frac{1}{3}\bigl(h_{ij}h_{k\ell}+h_{ik}h_{j\ell}+h_{i\ell}h_{jk}\bigr),\\
 X_{ijk\ell}
 &=\frac{S}{3}\bigl(h_{ij}K_{k\ell}+h_{ik}K_{j\ell}+h_{i\ell}K_{jk}
                   +K_{ij}h_{k\ell}+K_{ik}h_{j\ell}+K_{i\ell}h_{jk}\bigr).
\end{aligned}
\]
For every symmetric $4$-tensor $T$, we have $\langle T,Z\rangle=\langle T,h\otimes h\rangle$ and $\langle T,X\rangle=2S\langle T,h\otimes K\rangle$. Taking $T=u$ and substituting $E=h,K$ in \eqref{eq:weighted-trace-u}, we obtain
\begin{equation}\label{eq:u-zx-contractions}
 \langle u,Z\rangle=-\mathcal{C},\qquad
 \langle u,X\rangle=S^3|K|^2-2\left(\mathcal{A}-\frac{f_3\mathcal{C}}{S}-\frac{S}{n}|\nabla h|^2\right).
\end{equation}

In \eqref{eq:weighted-gradient}, $1-v$ is the normalized squared norm of $h(\nabla f_3)$. For $\nabla f_3$, set
\[
 b=
 \begin{cases}
 \displaystyle\left(1+\frac{2}{n}\right)\frac{|\nabla f_3|^2}{9S^2|PK|^2|\nabla h|^2}, & PK\ne0,\\[5pt]
 0, & PK=0.
 \end{cases}
\]

\begin{lem}\label{lem:gradient-contractions}
We have $0\le b\le v\le1$ and
\begin{equation}\label{eq:gradient-contraction-bounds}
 \frac{\mathcal C^2}{S|\nabla h|^4}+\frac{4nb^2|PK|^2}{3(n+2)(1+b-v)}
 \le\frac{\mathcal A+2\mathcal B}{3S|\nabla h|^2}
 \le\frac{\mathcal A}{S|\nabla h|^2}\le\frac{1}{2},
\end{equation}
where $\frac{b^2}{1+b-v}$ is assigned the value zero at $(b,v)=(0,1)$.
\end{lem}
\begin{proof}
For $PK\ne0$, \eqref{eq:component-projection-bounds} gives
\[
\begin{aligned}
 (1+b-v)|\nabla h|^2
 \le\sum_{i,k}\left(1+\frac{2}{n}+2\mu_k^2\right)h_{iik}^2\le\sum_k\left(h_{kkk}^2+3\sum_{i\ne k}h_{iik}^2\right)\le|\nabla h|^2,
\end{aligned}
\]
where the last sum counts the terms of $|\nabla h|^2$ with a repeated index. Thus $b\le v$. Away from $(0,1)$, the bound $0\le\frac{b^2}{1+b-v}\le b$ gives the stated extension by continuity.

When $\nabla f_3\ne0$, the Cauchy--Schwarz inequality gives
\[
 \sum_k\frac{(\nabla_k f_3)^2}{(n+2)S+2n\lambda_k^2}\ge\frac{|\nabla f_3|^4}{(n+2)S|\nabla f_3|^2+2n|h(\nabla f_3)|^2}.
\]
The first inequality in \eqref{eq:gradient-contraction-bounds} follows from \eqref{audit:lowerAB}, using this estimate when $\nabla f_3\ne0$ and $b=0$ otherwise.

The middle inequality follows from \eqref{eq:shape-commutator}.
If $\mu_i^2\le\frac{1}{2}$ for every $i$, then $\frac{\mathcal A}{S}=\sum_i\mu_i^2|\nabla_i h|^2\le\frac{1}{2}|\nabla h|^2$.
Otherwise choose $k$ with $\mu_k^2>\frac{1}{2}$, so that $\mu_i^2\le1-\mu_k^2$ for $i\ne k$. Then \eqref{audit:Tv} with $\eta=e_k$ gives $\frac{|\nabla_k h|^2}{|\nabla h|^2}\le\frac{1}{2}$. Hence
\[
 \frac{\mathcal A}{S|\nabla h|^2}
 \le 1-\mu_k^2+(2\mu_k^2-1)\frac{|\nabla_k h|^2}{|\nabla h|^2}
 \le\frac{1}{2} .
\]
\end{proof}

\section{Scalar inequalities and reduction to the endpoint}\label{sec:scalar-inequality}

Since $|\nabla F|(p_j)\to0$ along the maximizing sequence, \eqref{eq:trace-moment-L} gives $4(\mathcal A-\mathcal B)\le2S|\nabla h|^2$ for the limiting components. We therefore fix $L=2$ and combine this bound with the Laplacian inequality in \eqref{F} and the tensor estimates of Section~\ref{sec:pre} to obtain \eqref{master} and \eqref{dual}. With $\frac{h}{\sqrt S}$, $\frac{\nabla h}{|\nabla h|}$ and the weights fixed, monotonicity in the rescaling parameter $s$ reduces the gap estimate to $s=\frac{8n}{5}$. The same reduction holds when we retain the exact norm of $Z$ in the projection estimate.

For $\kappa>0$, set
\[
\begin{aligned}
 B_0&=\frac{4}{3}(3+2\kappa)-\frac{4\kappa}{n}-\frac{S(3-2\kappa)(1+2\kappa)}{2(S-n)},\\
 G&=\frac{S-n}{S}\left(2+\frac{3(n+4)}{8n^2}\right)-1+\frac{3-4\kappa}{3n}.
\end{aligned}
\]

\begin{lem}\label{lem:scalar-tests}
For the limiting components along $(p_j)$ and every $\kappa>0$, we have
\begin{equation}\label{master}
\begin{aligned}
 &\frac{\bigl(\mathcal C-\frac{S-n}{n}f_3+\frac{4\kappa S}{3}\tr(h^3K)\bigr)^2}{S^4(S-n)}
  +\frac{(n+2)(3+2\kappa)^2|h(\nabla f_3)|^2}{9S(n+2+3nv)|\nabla h|^4}\\
 +&\frac{8\kappa^2}{3(S-n)}\left(\frac{S}{2}|K|^2-\tr(h^2K^2)\right)
  +(15+4\kappa)\left(\frac{\mathcal C^2}{S|\nabla h|^4}
       +\frac{4nb^2|PK|^2}{3(n+2)(1+b-v)}\right)\\
 \le& G+B_0|K|^2+\frac{4(3+2\kappa)}{3}
       \left(\frac{f_3^2}{S^3}-\frac{\Delta(F-\frac{S}{4} f_4)}{S^2|\nabla h|^2}\right)
       -\frac{4\kappa f_3\mathcal C}{S^2|\nabla h|^2}.
\end{aligned}
\end{equation}
We also have
\begin{equation}\label{dual}
\begin{aligned}
  &|K|^2+\frac{3(n+2)(1-v)|PK|^2}{n+2+3nv}+\frac{2(S-n)}{S} \left(\frac{\mathcal{C}^2}{S|\nabla h|^4}+\frac{4nb^2|PK|^2}{3(n+2)(1+b-v)}\right)\\
  \le& \frac{2(S-n)}{3S}\left(1-\frac{3}{n}\right) +\frac{4(n+1)(n+3)}{3n(n+4)}\left(\frac{S-n}{S}\right)^2 .
\end{aligned}
\end{equation}
\end{lem}
\begin{proof}
Following~\cite{CLWY}, take $W=u-\kappa X+\frac{3(S-n)}{4nS}Y$. Substituting $E=h,K$ into \eqref{eq:weighted-trace-y} and \eqref{eq:weighted-trace-v} gives
\[
 \begin{aligned}
 \langle Y,Z\rangle=\frac{4S}{3}f_3,  &\qquad      \langle V,Z\rangle=0,\\
 \langle Y,X\rangle=\frac{8}{3}S^3|K|^2,   &\qquad \langle V,X\rangle=\frac{4|h(\nabla f_3)|^2}{9S^2}.
 \end{aligned}
\]
Replacing $h$ by $h+tK$ in the Gauss equation gives sectional curvature derivatives $\lambda_iK_{jj}+\lambda_jK_{ii}$ at $t=0$ for $i\ne j$. In the principal frame, $\langle h,K\rangle=0$ gives
\begin{align}
 2|K|^2-\frac{|X|^2}{S^3}
 &=\frac{8}{3S}\left(\frac{S}{2}|K|^2-\tr(h^2K^2)\right)
  =\frac{4}{3S}\sum_{i<j}(\lambda_iK_{jj}+\lambda_jK_{ii})^2, \label{xs}\\
 \langle X,Z\rangle&=\frac{4S}{3}\tr(h^3K). \label{zs}
\end{align}

Together with \eqref{eq:u-zx-contractions}, these identities give
\[
 \langle W,Z\rangle=-\mathcal{C}+\frac{S-n}{n}f_3-\frac{4\kappa S}{3}\tr(h^3K),
 \qquad \langle W,V\rangle=-\frac{2(3+2\kappa)|h(\nabla f_3)|^2}{9S^2}.
\]
When $h(\nabla f_3)\ne0$, these identities give $V\ne0$. Applying Bessel's inequality to $Z\perp V$ and using $|Z|^2=\frac{S^2+2f_4}{3}\le S^2$ and \eqref{eq:symmetric-test-norm}, we get
\[
 |W|^2\ge
 \frac{\bigl(\mathcal{C}-\frac{S-n}{n}f_3+\frac{4\kappa S}{3}\tr(h^3K)\bigr)^2}{S^2}
 +\frac{(n+2)(3+2\kappa)^2|h(\nabla f_3)|^2}{9(n+2+3nv)|\nabla h|^2}.
\]
If $h(\nabla f_3)=0$, then $V=0$ and projection onto $Z$ alone gives the same inequality.

At the limit $\Delta F\le0$ by \eqref{F}, so this lower bound also holds for $|W|^2-\frac{a\Delta F}{S}$ whenever $a>0$. On expanding, the coefficients of $\mathcal A$ and $\mathcal B$ are $3+4\kappa-2a$ and $-6-a$, respectively. Taking $a=\frac{4(3+2\kappa)}{3}$ makes the second twice the first. In $\Delta F$ we expand only $\frac{S}{4}\Delta f_4$ by \eqref{eq:trace-laplacians}. Using $f_4=S^2|K|^2+\frac{f_3^2}{S}+\frac{S^2}{n}$, \eqref{Uidentity} and \eqref{xs}, we obtain
\[
\begin{aligned}
 &\frac{|W|^2-\frac{a\Delta F}{S}}{S|\nabla h|^2}
  +\frac{15+4\kappa}{3}\,\frac{\mathcal A+2\mathcal B}{S|\nabla h|^2}
  +\frac{8\kappa^2}{3(S-n)}\left(\frac{S}{2}|K|^2-\tr(h^2K^2)\right)\\
 =&G+B_0|K|^2+a\left(\frac{f_3^2}{S^3}-\frac{\Delta(F-\frac{S}{4} f_4)}{S^2|\nabla h|^2}\right)
   -\frac{4\kappa f_3\mathcal C}{S^2|\nabla h|^2}.
\end{aligned}
\]
The projection bound for $|W|^2$ and \eqref{eq:gradient-contraction-bounds} now give \eqref{master}.

For \eqref{dual}, substitute $3(\mathcal A-2\mathcal B)=4(\mathcal A-\mathcal B)-(\mathcal A+2\mathcal B)$ in \eqref{Uidentity}. Using \eqref{eq:trace-moment-L} for $4(\mathcal A-\mathcal B)$ and \eqref{eq:gradient-contraction-bounds} for $\mathcal A+2\mathcal B$ gives an upper bound for $|u|^2$. Comparing it with \eqref{UVlower} and using \eqref{eq:simons-gradient} proves \eqref{dual}.
\end{proof}

Let $\mu_i$ and $T_{ijk}$ be the respective limits of $\frac{\lambda_i}{\sqrt{S}}$ and $\frac{h_{ijk}}{|\nabla h|}$ along $(p_j)$. For each $s>n$, set
\[
 \lambda_i(s)=\sqrt{s}\,\mu_i,\qquad
 h_{ijk}(s)=\sqrt{s(s-n)}\,T_{ijk},
\]
so that $\lambda_i(S)=\lambda_i$ and $h_{ijk}(S)=h_{ijk}$ by \eqref{eq:simons-gradient}.

Fix the weights $\kappa$ and $\lambda$ and add $\lambda$ times \eqref{dual} to \eqref{master}. For $s>n$, replace $S$ by $s$ in the two expressions and evaluate them using $\lambda_i(s)$ and $h_{ijk}(s)$. Let $\Phi(s)$ be the resulting right-hand side minus the left-hand side. We write $h$ and $\nabla h$ for the rescaled arrays and suppress the argument $s$ in the quantities computed from them.

\begin{lem}\label{lem:endpoint}
Let $n\ge4$, $0<\kappa\le\frac{3}{2}$ and $0\le\lambda\le\frac{3}{4}$. Then $\Phi$ is strictly increasing on $(n,\infty)$ and $\Phi(s)\ge\Phi(S)\ge0$ for every $s\ge S$.
\end{lem}
\begin{proof}
Since $\frac{f_3}{s^{3/2}}=\sum_i\mu_i^3$, the formulas for $F$ show that $\frac{\partial_{f_3}F}{s^{3/2}}$ and $\partial_{f_3}^2F$ depend only on the normalized principal curvatures. The diagonal entries of $K$ are $K_{ii}=\mu_i^2-\mu_i\sum_j\mu_j^3-\frac{1}{n}$, and rescaling the generators of $\ker P$ gives
\[
 \ker P=\operatorname{span}\left\{(1,\ldots,1),\,(\mu_i)_i,\,
 \left(\mu_i^3+\frac{3\,\partial_{f_3}F}{s^{3/2}}\mu_i^2\right)_i\right\}.
\]
Thus both $K$ and the orthogonal projection $P$ are independent of $s$.

The rescaling gives $|\nabla h|^2=s(s-n)$ and
\[
 \frac{\mathcal C}{\sqrt{s}\,|\nabla h|^2}=\sum_{i,j,k}\mu_i T_{ijk}^2,\qquad
 \frac{\nabla_k f_3}{3s|\nabla h|}=\sum_i\mu_i^2T_{iik}.
\]
The pair $(b,v)$ is independent of $s$, since $(b,v)=(0,1)$ for $PK=0$ and, for $PK\ne0$, both entries are determined by summing the squares of the second identity with weights $1$ and $\mu_k^2$. By \eqref{eq:trace-laplacians},
\[
 \frac{\Delta f_3}{3\sqrt{s}\,|\nabla h|^2}
 =2\sum_{i,j,k}\mu_i T_{ijk}^2-\sum_i\mu_i^3.
\]
The chain rule therefore shows that $\frac{\Delta(F-\frac{s}{4}f_4)}{s^2|\nabla h|^2}$ is independent of $s$.

The contribution to $\Phi(s)$ from the projection onto $Z$ in \eqref{master} is
\[
 -\frac{s}{s-n}\left[
 \frac{s-n}{s}\left(\sum_{i,j,k}\mu_i T_{ijk}^2-\frac{1}{n}\sum_i\mu_i^3\right)
 +\frac{4\kappa}{3}\sum_i\mu_i^3K_{ii}\right]^2 .
\]
Expanding this square and collecting powers of $\frac{s-n}{s}$ in $\Phi(s)$ gives coefficients depending only on $n$, the weights $\kappa,\lambda$ and the normalized components $\mu_i,T_{ijk}$. The expansion and its derivative are
\[
 \begin{aligned}
 \Phi(s)&=\sigma_{-1}\frac{s}{s-n}+\sigma_0+\sigma_1\frac{s-n}{s}
          +\sigma_2\left(\frac{s-n}{s}\right)^2,\\
 \Phi'(s)&=-\frac{n\sigma_{-1}}{(s-n)^2}+\frac{n\sigma_1}{s^2}
           +\frac{2n(s-n)\sigma_2}{s^3},
 \end{aligned}
\]
where $\sigma_2=\frac{4\lambda(n+1)(n+3)}{3n(n+4)}\ge0$ by \eqref{dual}. Substituting \eqref{eq:weighted-gradient} gives
\[
\begin{aligned}
 -\sigma_{-1}={}&\left(\frac{4\kappa}{3}\sum_i\mu_i^3K_{ii}\right)^2
 +\frac{(3-2\kappa)(1+2\kappa)}{2}|K|^2\\
 &+\frac{(n+2)(3+2\kappa)^2(1-v)}{2(n+2+3nv)}|PK|^2
 +\frac{8\kappa^2}{3}\left(\frac{|K|^2}{2}-\sum_i\mu_i^2K_{ii}^2\right)\ge0,
\end{aligned}
\]
where the last term is nonnegative by \eqref{xs}. The linear coefficient is
\[
\begin{aligned}
 \sigma_1={}&2+\frac{3(n+4)}{8n^2}+\frac{2\lambda}{3}\left(1-\frac{3}{n}\right)
 -\left(\sum_{i,j,k}\mu_i T_{ijk}^2-\frac{1}{n}\sum_i\mu_i^3\right)^2\\
 &-2\lambda\left[\left(\sum_{i,j,k}\mu_i T_{ijk}^2\right)^2
 +\frac{4nb^2|PK|^2}{3(n+2)(1+b-v)}\right].
\end{aligned}
\]
Passing to the limit in \eqref{eq:gradient-contraction-bounds} along $(p_j)$ shows that the bracket multiplying $-2\lambda$ is at most $\frac{1}{2}$, so $\left|\sum_{i,j,k}\mu_i T_{ijk}^2\right|\le\frac{1}{\sqrt{2}}$. Using $\left|\sum_i\mu_i^3\right|\le1$ and $\lambda(1-\frac{3}{n})\ge0$, we obtain
\[
 \sigma_1\ge 2+\frac{3(n+4)}{8n^2}
 -\left(\frac{1}{\sqrt{2}}+\frac{1}{n}\right)^2-\lambda
 >1-\lambda\ge\frac{1}{4}.
\]
Thus $\Phi'(s)>0$ for $s>n$. At $s=S$, adding $\lambda$ times \eqref{dual} to \eqref{master} gives $\Phi(S)\ge0$. Monotonicity then yields $\Phi(s)\ge\Phi(S)\ge0$ for all $s\ge S$.
\end{proof}

The $Z$-projection in Bessel's inequality has denominator $|Z|^2$, and any larger denominator also gives a lower bound. For the same weights and the rescaled data, retaining the exact denominator multiplies the $Z$-projection term in \eqref{master} by $\frac{s^2}{|Z|^2}$. The resulting weighted difference is
\[
 \Phi_Z(s)=\Phi(s)-\left(\frac{s^2}{|Z|^2}-1\right)
 \frac{\bigl(\mathcal C-\frac{s-n}{n}f_3+\frac{4\kappa s}{3}\tr(h^3K)\bigr)^2}{s^4(s-n)}.
\]
At $s=S$, adding $\lambda$ times \eqref{dual} to this version of \eqref{master} gives $\Phi_Z(S)\ge0$.

The factor $\frac{s^2}{|Z|^2}=\frac{3}{1+2\sum_i\mu_i^4}\ge1$ is independent of $s$. Thus $\Phi_Z$ also expands in the powers $-1,0,1,2$ of $\frac{s-n}{s}$, with coefficients $\sigma_j^Z$. This factor multiplies the projection square in $-\sigma_{-1}$, and the quadratic term comes entirely from \eqref{dual}. Hence $\sigma_{-1}^Z\le\sigma_{-1}\le0$ and $\sigma_2^Z=\sigma_2\ge0$.
Multiplying the first square in $\sigma_1$ by $\frac{s^2}{|Z|^2}$ gives $\sigma_1^Z$.

At $s=S$, \eqref{eq:gradient-contraction-bounds} bounds the bracket multiplying $-2\lambda$ by $\frac{1}{2}$ and gives $|\mathcal C|\le\sqrt{\frac{S}{2}}\,|\nabla h|^2$. Since $f_3\ge0$ and $\frac{3|Z|^2}{S^2}=1+\frac{2}{n}+\frac{2f_3^2}{S^3}+2|K|^2$, applying the triangle inequality and then the Cauchy--Schwarz inequality gives
\[
 \frac{S}{|Z|^2}\left(\frac{\mathcal C}{|\nabla h|^2}-\frac{f_3}{nS}\right)^2
 \le\frac{3\left(\frac{1}{\sqrt 2}+\frac{f_3}{nS^{3/2}}\right)^2}
          {1+\frac{2}{n}+\frac{2f_3^2}{S^3}}
 \le\frac{3n}{2(n+2)}+\frac{3}{2n^2}.
\]
Substituting these bounds in $\sigma_1^Z$ and using $n\ge4$ and $\lambda\le\frac{3}{4}$, we obtain
\[
 \sigma_1^Z\ge\frac{1}{2}-\frac{\lambda}{3}
 +\frac{1}{n}\left(\frac{3n}{n+2}+\frac{3}{8}-2\lambda\right)
 \ge\frac{1}{4}.
\]
The derivative formula now gives $\Phi_Z'(s)>0$ for $s>n$. If $n<S\le\frac{8n}{5}$, monotonicity and the bound at $s=S$ give
\begin{equation}\label{eq:exact-endpoint-sign}
 \Phi_Z\!\left(\frac{8n}{5}\right)\ge\Phi_Z(S)\ge0 .
\end{equation}

To bound $\Phi_Z(\frac{8n}{5})$ from above, we first estimate the terms involving $b$ or $v$ uniformly over $0\le b\le v\le1$. Expanding $\Delta(F-\frac{s}{4}f_4)$ by the chain rule and substituting the definitions of $b,v$, we collect these terms in $\Phi(s)$ and in $\Phi_Z(s)$. In each case their sum is the same multiple of $|PK|^2$. The following lemma bounds this sum by $\rho^2|PK|^2$.
For $n\ge4$, $s>n$, $\kappa>0$ and $\lambda\ge0$, set
\[
 Q=\frac{4}{3}\left(15+4\kappa+\frac{2\lambda(s-n)}{s}\right),\qquad
 C_\lambda=\frac{s(3+2\kappa)^2}{s-n}+6\lambda.
\]
When $0<-12(3+2\kappa)\partial_{f_3}^2F<Q$, define
\[
 \rho=\frac{1}{\sqrt{6}}\left[2\sqrt{\frac{2n+1}{n+2}}\,\bigl(\sqrt Q-\sqrt{Q+12(3+2\kappa)\partial_{f_3}^2F}\bigr)-\sqrt{\frac{(n+2)C_\lambda}{n}}\right]_+ .
\]

\begin{lem}\label{lem:optimized-remainder}
For these parameters and all $0\le b\le v\le1$, the common coefficient satisfies
\begin{equation}\label{eq:optimized-aux}
 \frac{n}{n+2}\left(-12(3+2\kappa)b\,\partial_{f_3}^2F-\frac{Qb^2}{1+b-v}\right)-\frac{(n+2)C_\lambda(1-v)}{2(n+2+3nv)}\le \rho^2 ,
\end{equation}
where $\frac{b^2}{1+b-v}$ takes the value zero at $(b,v)=(0,1)$.
\end{lem}
\begin{proof}
At $(b,v)=(0,1)$ the left-hand side vanishes. Otherwise, completing the square in $b$ gives
\[
 -12(3+2\kappa)b\,\partial_{f_3}^2F-\frac{Qb^2}{1+b-v}
 \le (1-v)\bigl(\sqrt Q-\sqrt{Q+12(3+2\kappa)\partial_{f_3}^2F}\bigr)^2 .
\]
Indeed, multiplying the difference between the right- and left-hand sides by $1+b-v$ gives
\[
 \bigl((1+b-v)\sqrt{Q+12(3+2\kappa)\partial_{f_3}^2F}-(1-v)\sqrt Q\bigr)^2 .
\]
The coefficient of $|PK|^2$ is therefore at most
\[
 (1-v)\left[\frac{n}{n+2}\bigl(\sqrt Q-\sqrt{Q+12(3+2\kappa)\partial_{f_3}^2F}\bigr)^2
 -\frac{(n+2)C_\lambda}{2(n+2+3nv)}\right].
\]
For $a,c\ge0$ and $0\le v\le1$,
\[
 (1-v)\left(a-\frac{c}{n+2+3nv}\right)
 \le \frac{1}{3n}\left[\sqrt{(4n+2)a}-\sqrt c\right]_+^2 .
\]
If $(4n+2)a\le c$, the left-hand side is nonpositive. Otherwise, $3n$ times the difference between the right-hand and left-hand sides is
\[
 \left(\sqrt{(n+2+3nv)a}-\sqrt{\frac{(4n+2)c}{n+2+3nv}}\right)^2 .
\]
Taking $a=\frac{n}{n+2}(\sqrt Q-\sqrt{Q+12(3+2\kappa)\partial_{f_3}^2F})^2$ and $c=\frac{n+2}{2}C_\lambda$ gives \eqref{eq:optimized-aux}.
\end{proof}

\section{The spectral reduction and proof of the main theorem}\label{sec:spectral-reduction}

We use the principal curvatures to choose weights for which $\Phi_Z(\frac{8n}{5})$ has a strictly negative upper bound, contradicting \eqref{eq:exact-endpoint-sign}.
To estimate the terms in $|PK|^2$ and $\tr(h^2K^2)$ in \eqref{master}, consider multiplication by $\frac{h}{\sqrt{s}}$ on the space of diagonal operators. Its eigenvalues $\mu_i$ give trace zero, squared Hilbert--Schmidt norm one and cubic trace $\frac{f_3}{s^{3/2}}$.
In an orthonormal basis beginning with $\frac{I}{\sqrt{n}}$ and $\frac{h}{\sqrt{s}}$, the first two columns of its matrix are determined by
\[
 \frac{h}{\sqrt{s}}\,\frac{I}{\sqrt{n}}=\frac{h}{\sqrt{ns}},\qquad
 \frac{h^2}{s}=\frac{I}{n}+\frac{f_3}{s^2}h+K .
\]
Since this matrix is symmetric and $K$ is orthogonal to both $I$ and $h$, the compression to $\{I,h\}^{\perp}$ has trace $-\frac{f_3}{s^{3/2}}$ and squared Hilbert--Schmidt norm $1-\frac{2}{n}-\frac{f_3^2}{s^3}-2|K|^2$. The Cauchy--Schwarz inequality in dimension $n-2$ gives
\begin{equation}\label{eq:admissible-spectrum}
 2|K|^2+\frac{n-1}{n-2}\frac{f_3^2}{s^3}\le1-\frac{2}{n},\qquad \frac{f_3^2}{s^3}\le\frac{(n-2)^2}{n(n-1)} .
\end{equation}
When $K\ne0$, we have
\begin{equation}\label{eq:spectral-parameters}
 \frac{\tr(h^3K)}{s^{3/2}|K|^2}
 =\frac{f_3}{s^{3/2}}+\frac{\tr(hK^2)}{\sqrt{s}\,|K|^2}
 \ge\frac{3f_3}{2s^{3/2}}-1,
\end{equation}
where the equality follows by multiplying the formula for $\frac{h^2}{s}$ by $h$ and taking the inner product with $K$. For $\mu_i=\min_j\mu_j<0$, the quotient $\frac{\tr(hK^2)}{\sqrt{s}\,|K|^2}$ is at least $\mu_i$, since it is a weighted mean of the $\mu_j$ with weights $\frac{K_{jj}^2}{|K|^2}$. Using $|\mu_j|\le1$, we obtain
\[
 \frac{f_3}{s^{3/2}}\le\sum_{j\ne i}|\mu_j|^3
 \le1-\mu_i^2\le2(1+\mu_i).
\]
This proves the lower bound in \eqref{eq:spectral-parameters}. The contribution of the curvature term in \eqref{master} to $\Phi_Z(s)$ is
\[
 -\frac{8\kappa^2}{3(s-n)}\left(\frac{s}{2}|K|^2-\tr(h^2K^2)\right)
 =-\frac{8s\kappa^2}{3(s-n)}\left(\frac{|K|^2}{2}-\frac{|hK|^2}{s}\right).
\]

\begin{lem}\label{lem:spectral-coupling}
Let $h=\sqrt s\,\operatorname{diag}(\mu_1,\ldots,\mu_n)$, where $s>n\ge4$, $\sum_i\mu_i=0$ and $\sum_i\mu_i^2=1$. If $K\ne0$, then for all $a,c\ge0$,
\begin{equation}\label{eq:joint-spectral-bound}
\begin{aligned}
 a^2|PK|^2-c^2\left(\frac{|K|^2}{2}-\frac{|hK|^2}{s}\right)
 \le& |K|^2\left(a-\frac{c\,\bigl|\tr(h^3K)+3|K|^2\,\partial_{f_3}F\bigr|}{s^{3/2}|K|^2}\right)_+^2\\
 &-c^2|K|^2\left(|K|^2+\frac{f_3^2}{s^3}+\frac{2}{n}-\frac{1}{2}\right).
\end{aligned}
\end{equation}
\end{lem}
\begin{proof}
Let $\pi\colon\mathbb R^n\to\{I,h,K\}^{\perp}$ be the orthogonal projection. Since $hK\perp I$ and $\langle hK,h\rangle=s|K|^2$, the orthogonal decomposition is
\[
 hK=|K|^2h
 +\frac{\tr(hK^2)}{|K|^2}K+\pi(hK) .
\]
In an orthonormal basis of $\{I,h\}^{\perp}$ beginning with $\frac{K}{|K|}$, the compression of multiplication by $\frac{h}{\sqrt s}$ has the coordinates of $\frac{\pi(hK)}{\sqrt{s}\,|K|}$ in the off-diagonal positions of its first row and column. Its squared Hilbert--Schmidt norm satisfies
\[
 1-\frac{2}{n}-\frac{f_3^2}{s^3}-2|K|^2
 =\frac{\bigl(\tr(hK^2)\bigr)^2}{s|K|^4}
 +\frac{2|\pi(hK)|^2}{s|K|^2}+\frac{\|\pi h\pi\|_{\mathrm{HS}}^2}{s}.
\]
Taking squared norms in the decomposition and eliminating the squared trace term gives
\[
\begin{aligned}
 \frac{|K|^2}{2}-\frac{|hK|^2}{s}
 &=\frac{|\pi(hK)|^2}{s}+|K|^2\left(|K|^2+\frac{f_3^2}{s^3}+\frac{2}{n}-\frac{1}{2}\right)
   +\frac{|K|^2}{s}\|\pi h\pi\|_{\mathrm{HS}}^2\\
 &\ge\frac{|\pi(hK)|^2}{s}+|K|^2\left(|K|^2+\frac{f_3^2}{s^3}+\frac{2}{n}-\frac{1}{2}\right).
\end{aligned}
\]

Substituting the decompositions of $\frac{h^2}{s}$ and $hK$ into $P\bigl(h^3+\frac{3\,\partial_{f_3}F}{s}h^2\bigr)=0$ gives
\[
 \frac{\tr(h^3K)+3|K|^2\,\partial_{f_3}F}{s^{3/2}|K|^2}\,PK+\frac{P(\pi(hK))}{\sqrt{s}}=0 .
\]
For $|\tau|\le c$, the contraction property of $P$ and $K\perp\pi(hK)$ then give
\[
 a^2|PK|^2\le
 \left[a-\tau\,\frac{\tr(h^3K)+3|K|^2\,\partial_{f_3}F}{s^{3/2}|K|^2}\right]^2|K|^2+\frac{c^2|\pi(hK)|^2}{s} .
\]
Minimizing over $|\tau|\le c$ and using the preceding bound for $\frac{|\pi(hK)|^2}{s}$ proves \eqref{eq:joint-spectral-bound}.
\end{proof}

The sum of the $Z$-projection and curvature terms in $\Phi_Z(s)$ satisfies
\begin{equation}\label{eq:dimension-gain}
\begin{aligned}
 &-\frac{s-n}{|Z|^2}\left(\frac{\mathcal C}{|\nabla h|^2}-\frac{f_3}{ns}
      +\frac{4\kappa\,\tr(h^3K)}{3(s-n)}\right)^2
 -\frac{8s\kappa^2}{3(s-n)}\left(\frac{|K|^2}{2}-\frac{|hK|^2}{s}\right)\\
 \le& -\frac{s-n}{s^2}\left(\frac{\mathcal C}{|\nabla h|^2}-\frac{f_3}{ns}\right)^2
 -\frac{8s\kappa^2|K|^2}{3n(s-n)}.
\end{aligned}
\end{equation}
For $K\ne0$, the vectors $(\mu_i)_i$ and $\bigl(\frac{K_{ii}}{|K|}\bigr)_i$ are orthonormal, and the components of each sum to zero. Thus for $q\in\mathbb R$ the symmetric matrix with entries
\[
 B_{ij}=\frac{\mu_iK_{jj}+\mu_jK_{ii}}{|K|}-q\mu_i\mu_j
\]
has zero row sums, and at $q=0$ its off-diagonal entries are the normalized sectional curvature variations in \eqref{xs}. Multiplying the $i$th row sum by $B_{ii}-\frac{\tr B}{n}$ and summing over $i$ gives
\[
 \sum_i\left(B_{ii}-\frac{\tr B}{n}\right)^2
 =-\sum_{i<j}\left(B_{ii}+B_{jj}-\frac{2\tr B}{n}\right)B_{ij}.
\]
The squared coefficients of $B_{ij}$ on the right sum to $(n-2)\sum_i(B_{ii}-\frac{\tr B}{n})^2$. The Cauchy--Schwarz inequality therefore gives
\[
 \sum_i B_{ii}^2-\frac{(\tr B)^2}{n}\le(n-2)\sum_{i<j}B_{ij}^2 .
\]
Adding the off-diagonal squares and using $\|B\|^2=2+q^2$ and $\tr B=-q$, we obtain
\[
 \sum_{i\ne j} B_{ij}^2\ge \frac{2}{n}\left(\|B\|^2-\frac{(\tr B)^2}{n}\right)\ge\frac{4}{n} .
\]
Expanding $\sum_{i\ne j}B_{ij}^2$ gives
\[
 4\left(\frac{1}{2}-\frac{|hK|^2}{s|K|^2}\right)
 +\frac{4q\,\tr(h^3K)}{s^{3/2}|K|}+q^2\left(1-\frac{f_4}{s^2}\right)\ge\frac{4}{n} .
\]
At least two $\mu_i$ are nonzero, so $f_4<s^2$. Minimizing in $q$ gives
\[
 \frac{1}{2}-\frac{|hK|^2}{s|K|^2}\ge\frac{1}{n}+\frac{\bigl(\tr(h^3K)\bigr)^2}{s|K|^2(s^2-f_4)}.
\]
Since $s^2-|Z|^2=\frac{2}{3}(s^2-f_4)$, completing the square bounds the $Z$-projection term in \eqref{eq:dimension-gain} from above by
\[
 -\frac{s-n}{s^2}\left(\frac{\mathcal C}{|\nabla h|^2}-\frac{f_3}{ns}\right)^2
 +\frac{8\kappa^2\bigl(\tr(h^3K)\bigr)^2}{3(s-n)(s^2-f_4)} .
\]
The squared trace terms cancel on adding the curvature estimate, giving \eqref{eq:dimension-gain}. For $K=0$ it follows directly from $|Z|^2\le s^2$.

To apply Lemma~\ref{lem:optimized-remainder} at $s=\frac{8n}{5}$, set
\[
 B=Q+12(3+2\kappa)\partial_{f_3}^2F
   =12-\left(\frac{211}{75}-\frac{6f_3^2}{5s^3}\right)(3+2\kappa)+\lambda,
\]
where $B<Q$ since $f_3^2\le s^3$. For $B>0$ the lemma bounds the $b,v$ terms by $\rho^2|PK|^2$.

Since $s-n=\frac{3s}{8}$, the coefficients in the projection and curvature terms of \eqref{eq:dimension-gain} satisfy
\[
-\frac{s-n}{|Z|^2}=-\frac{3s}{8|Z|^2},\qquad
 \frac{4\kappa}{3(s-n)}=\frac{32\kappa}{9s},\qquad
 -\frac{8s\kappa^2}{3(s-n)}=-\frac{64\kappa^2}{9}.
\]
\begingroup
\setlength{\emergencystretch}{\linewidth}
Expanding the Laplacian term in \eqref{master} by the chain rule
and applying the $b,v$ estimate gives $\mathcal E$ as an upper bound
for the sum of the remaining terms of $\Phi_Z(\frac{8n}{5})$. Thus
\begin{subequations}
\begin{equation}\label{eq:endpoint-upper-bound}
\begin{aligned}
 \Phi_Z\!\left(\frac{8n}{5}\right)
 \le \mathcal E-\frac{3s}{8|Z|^2}
 \left(\frac{\mathcal C}{|\nabla h|^2}-\frac{f_3}{ns}
       +\frac{32\kappa\,\tr(h^3K)}{9s}\right)^2-\frac{64\kappa^2}{9}
 \left(\frac{|K|^2}{2}-\frac{|hK|^2}{s}\right),
\end{aligned}
\end{equation}
where
\begin{align}
 \mathcal E={}&G+\frac{4(3+2\kappa)}{3}\left(\frac{f_3^2}{s^3}
 -\frac{(\partial_{f_3}F)\Delta f_3}{s^2|\nabla h|^2}\right)
 -\frac{4\kappa f_3\mathcal C}{s^2|\nabla h|^2}\notag\\
 &-\frac{3Q\mathcal C^2}{4s|\nabla h|^4}
 +(B_0-\lambda)|K|^2+\rho^2|PK|^2+\lambda\left(\frac{7}{16}-\frac{39}{64n}-\frac{9}{64(n+4)}\right),\label{eq:endpoint-common}\\
 G={}&-\frac{1}{4}+\frac{1}{n}\left(\frac{73}{64}-\frac{4\kappa}{3}\right)
      +\frac{9}{16n^2},\qquad
 B_0=\frac{16\kappa}{3}\left(\kappa-\frac{1}{2}-\frac{3}{4n}\right),\notag
\end{align}
\end{subequations}
with $\Delta f_3$ given by \eqref{eq:trace-laplacians}.

\par\endgroup

\begin{lem}\label{lem:endpoint-verification}
Let $n\ge4$ and fix principal curvatures $\lambda_i$ with $\sum_i\lambda_i=0$, $\sum_i\lambda_i^2=s=\frac{8n}{5}$ and $f_3=\sum_i\lambda_i^3\ge0$.
There exist weights \mbox{$0<\kappa\le\frac{3}{2}$} and \mbox{$0\le\lambda\le\frac{11}{20}$}, depending only on $n$ and $\frac{f_3^2}{s^3}$, such that $B>0$. For these weights and every value of $\mathcal C=\sum_{i,j,k}\lambda_i h_{ijk}^2$, we have
\[
\begin{aligned}
 \mathcal E-\frac{3s}{8|Z|^2}
 \left(\frac{\mathcal C}{|\nabla h|^2}-\frac{f_3}{ns}
       +\frac{32\kappa\,\tr(h^3K)}{9s}\right)^2-\frac{64\kappa^2}{9}
 \left(\frac{|K|^2}{2}-\frac{|hK|^2}{s}\right)<-10^{-4},
\end{aligned}
\]
where $\mathcal E$ is given by \eqref{eq:endpoint-common}.
\end{lem}

\begin{proof}[Proof of Lemma~\ref{lem:endpoint-verification}]
For $\frac{f_3^2}{s^3}\le\frac{2}{5}$, choose
\[
 \kappa=\frac{1}{15}+\frac{2f_3^2}{5s^3}
       +\frac{4}{3n}\left(1+\frac{5f_3^2}{2s^3}\right)\in(0,1),\quad
 \lambda=\left(\frac{11}{20}-\frac{51f_3^2}{50s^3}\right)\left(1-\frac{4}{n}\right)
         \in\left[0,\frac{11}{20}\right].
\]
Applying \eqref{eq:dimension-gain} to the projection and curvature terms in \eqref{eq:endpoint-upper-bound} and completing the square in $\mathcal C$ gives the upper bound
\begin{gather}
 \mathcal P-T|K|^2+\rho^2|PK|^2,\notag\\
 \mathcal P=G+\frac{f_3^2}{s^3}\left(J-\frac{3}{8n^2}\right)+\lambda\left(\frac{7}{16}-\frac{39}{64n}-\frac{9}{64(n+4)}\right),\label{eq:small-bound}\displaybreak[0]\\
 T=\lambda+\frac{16\kappa}{3}
 \left[\frac{1}{2}+\frac{3}{4n}-\left(1-\frac{4}{3n}\right)\kappa\right],\notag\\
 J=\max_{p\in\mathbb R}\left\{\frac{4(3+2\kappa)}{3}
 \left[1+(2p-1)\left(\frac{137}{100}-\frac{f_3^2}{10s^3}\right)\right]
 +2p\left(\frac{3}{8n}-2\kappa\right)-\left(\frac{3Q}{4}+\frac{3}{8}\right)p^2\right\}.\notag
\end{gather}
Since $|PK|\le |K|$, it suffices to prove $\mathcal P<-10^{-4}$, $B>0$ and $\rho^2<T$. We first bound $\mathcal P$.

At $p=\frac{1}{2}$ the derivative is at most $12+4\kappa-\frac{3Q}{4}\le-3$. Hence the maximum lies at $p\le\frac{1}{2}$. For each such $p$, the bracket is affine in $\frac{1}{n}$ and its quadratic coefficient in $\frac{f_3^2}{s^3}$ is $\frac{4}{15}(1-2p)(\frac{2}{5}+\frac{10}{3n})\ge0$.
At the four corners $(\frac{1}{n},\frac{f_3^2}{s^3})=(0,0),(\frac{1}{4},0),(0,\frac{2}{5}),(\frac{1}{4},\frac{2}{5})$,
the maximum in $p$ is at most $0.401$, $0.416$, $0.443$ and $0.332$, respectively.
Taking chords in the two variables gives
\[
 J\le0.401+\frac{0.06}{n}+0.105\left(1-\frac{12}{n}\right)\frac{f_3^2}{s^3},
\]
with all decimals in the proof understood as exact.

\begingroup
\setlength{\emergencystretch}{\linewidth}
Evaluating the right-hand side of \eqref{dual} at $\frac{s-n}{s}=\frac{3}{8}$
gives the coefficient of $\lambda$ in $\mathcal P$. Using $n+4\le2n$
to bound its last term from above, we obtain
\[
\frac{7}{16}-\frac{39}{64n}-\frac{9}{64(n+4)}
 \le 0.4375-\frac{0.6796}{n}.
\]
\par\endgroup
The resulting quadratic upper bound for $\mathcal P$ has leading coefficient $0.105(1-\frac{12}{n})$ in $\frac{f_3^2}{s^3}$. Its values at $0,\frac{2}{5}$ are at most $-0.009-\frac{0.16}{n}$ and $-0.01+\frac{0.32}{n}-\frac{2.72}{n^2}$. Since $n\ge4$, its excess over the chord is at most $\frac{0.1344}{n}(1-\frac{5f_3^2}{2s^3})$, which is absorbed by the first endpoint term. Hence
\[
 \mathcal P\le\max\left\{-0.009,\,-0.01+\frac{0.32}{n}-\frac{2.72}{n^2}\right\}
 \le-\frac{1}{1700}.
\]

To prove $\rho^2<T$, use $\sqrt{(8+\frac{4}{n})Q}\le\frac{4Q}{13}+\frac{13}{4}(2+\frac{1}{n})$. By Lemma~\ref{lem:optimized-remainder}, once $B,T>0$ it suffices to show
\[
 \sqrt{\left(8+\frac{4}{n}\right)B}
 +\left(1+\frac{2}{n}\right)\sqrt{C_\lambda}
 +\sqrt{6\left(1+\frac{2}{n}\right)T}
 >\frac{4Q}{13}+\frac{13}{4}\left(2+\frac{1}{n}\right).
\]
In the following interpolation, subscripts $-,+$ denote the limit as $n\to\infty$ and the value at $n=4$, with $\frac{f_3^2}{s^3}$ fixed. As this ratio increases from $0$ to $\frac{2}{5}$, $8B_-,C_-$ increase and $6T_-$ decreases. Their endpoint values give $5.4\le\sqrt{8B_-},\sqrt{C_-}\le5.8$ and $1.6\le\sqrt{6T_-}\le2.1$. Applying the secant estimates to the three radicals at $n\to\infty$ gives
\[
 \begin{aligned}
 &\sqrt{8B_-}+\sqrt{C_-}+\sqrt{6T_-}-\frac{4Q_-}{13}-\frac{13}{2}\\
\ge& \frac{8B_-+C_-+62.64}{11.2}+\frac{6T_-+3.36}{3.7}
 -\frac{4Q_-}{13}-\frac{13}{2}>\frac{1}{200},
 \end{aligned}
\]
where the last inequality follows because the secant lower bound is concave in $\frac{f_3^2}{s^3}$ and exceeds $\frac{1}{200}$ at both endpoints.

Both $B$ and $Q$ are affine in $\frac{1}{n}$, and $B_+$ decreases with $\frac{f_3^2}{s^3}$ to $\frac{187}{225}>0.9^2$. Since $4(\kappa_+-\kappa_-)\le\frac{8}{3}<3+2\kappa_-$, we have $3+2\kappa\le(1+\frac{2}{n})(3+2\kappa_-)$ and $\sqrt{C_\lambda}\le(1+\frac{2}{n})\sqrt{C_-}$. Rationalizing $\sqrt{C_\lambda}-\sqrt{\frac{8}{3}}(3+2\kappa)$ and using $\lambda=(1-\frac{4}{n})\lambda_-$ gives
\[
 \left(1+\frac{2}{n}\right)\sqrt{C_\lambda}
 \ge\left(1-\frac{4}{n}\right)\sqrt{C_-}
 +\frac{4}{n}\sqrt{\frac{8}{3}}\left(\frac{9}{2}+\kappa_-+2\kappa_+\right).
\]
For $T$, subtraction of its endpoint chord gives
\[
 \begin{aligned}
 T-\left(1-\frac{4}{n}\right)T_--\frac{4}{n}T_+=\frac{64}{9n}(\kappa_+-\kappa)
 \left[(\kappa_+-\kappa)-\left(3\kappa_-+\frac{9}{16}-\kappa_+\right)\right]\ge-\frac{1}{4n},
 \end{aligned}
\]
where the lower bound follows by completing the square in the nonnegative difference $\kappa_+-\kappa$, using $3\kappa_-+\frac{9}{16}-\kappa_+=\frac{29}{80}-\frac{f_3^2}{30s^3}<\frac{3}{8}$. Since $T_+=\kappa_+(\frac{11}{3}-\frac{32\kappa_+}{9})$, the concave quadratic $T_+-(0.9-\frac{5f_3^2}{6s^3})^2$ is at least $\frac{79}{900}>\frac{1}{16}$ at both endpoints of the trace interval. Thus
\[
 T\ge\left(1-\frac{4}{n}\right)T_-+\frac{4}{n}
 \left(0.9-\frac{5f_3^2}{6s^3}\right)^2>0 .
\]
By concavity of the geometric mean for the first and third radicals and the lower bound in the limit as $n\to\infty$, the sum of the three radicals exceeds the affine upper bound for $\sqrt{(8+\frac{4}{n})Q}$ by at least
\[
 \frac{1}{200}\left(1-\frac{4}{n}\right)+
 \frac{4}{n}\left[5.4-\frac{5f_3^2}{2s^3}
 +1.63\left(\frac{9}{2}+\kappa_-+2\kappa_+\right)-\frac{4Q_+}{13}-\frac{117}{16}\right]>0,
\]
where $1.63<\sqrt{\frac{8}{3}}$ and the bracket is affine in $\frac{f_3^2}{s^3}$ with both endpoint values greater than $\frac{1}{50}$. Since $\rho^2<T$ and $|PK|\le|K|$, we have
\[
 \mathcal P-T|K|^2+\rho^2|PK|^2\le\mathcal P<-10^{-4}.
\]

For $\frac{f_3^2}{s^3}\ge\frac{2}{5}$ take $\kappa=\frac{3f_3^2}{2s^3}\in[\frac{3}{5},\frac{3}{2}]$ and $\lambda=0$. In this case $n\ge5$, since \eqref{eq:admissible-spectrum} gives $\frac{f_3^2}{s^3}\le\frac{1}{3}$ for $n=4$. For $Q=20+\frac{16\kappa}{3}$ and $B=\frac{89}{25}-\frac{242\kappa}{75}+\frac{8\kappa^2}{5}$, we have
\[
 \sqrt Q\le4.82+\frac{5}{9}\left(\kappa-\frac{3}{5}\right),\qquad
 \sqrt B\ge1.67-\frac{\kappa}{3}>0 .
\]
After squaring, both bounds follow from quadratics with positive leading coefficient and negative discriminant. By Lemma~\ref{lem:optimized-remainder},
\begin{equation}\label{eq:root-bracket}
 \frac{\rho}{\sqrt\kappa}\le\frac{1}{\sqrt{6\kappa}}\left(\frac{37}{12}-\frac{13\kappa}{18}-\frac{6}{n}\right).
\end{equation}
Indeed $\sqrt{8}<2.83$, $\sqrt{1+\frac{2}{n}}\ge1+\frac{0.9}{n}$ and $\sqrt{\frac{8}{3}}>1.63$ bound the numerator before taking the positive part by $2.83[3.35+\frac{8}{9}(\kappa-\frac{3}{5})]-1.63(3+2\kappa)(1+\frac{0.9}{n})$.

Take $q_0=\frac{3}{2}-\frac{5\kappa}{6}$ and $q=(1-\frac{4}{n})q_0$. The tangent line bounds the $Z$-projection term in \eqref{eq:endpoint-upper-bound} from above by
\[
 -\frac{3qf_3}{4s^2}\left(\frac{\mathcal C}{|\nabla h|^2}-\frac{f_3}{ns}
            +\frac{32\kappa\,\tr(h^3K)}{9s}\right)
 +\frac{3q^2f_3^2|Z|^2}{8s^5} .
\]
Using $\frac{3|Z|^2}{s^2}=1+\frac{2}{n}+\frac{4\kappa}{3}+2|K|^2$, complete the square in $\frac{\mathcal C}{|\nabla h|^2}-\frac{f_3}{3s}$. For $K=0$ the upper bound is $\mathcal P$. For $K\ne0$, Lemma~\ref{lem:spectral-coupling} with $a=\rho$ and $c=\frac{8\kappa}{3}$ then gives the upper bound
\[
 \mathcal P+\kappa|K|^2\left[
 V\!\left(\frac{\tr(h^3K)}{s^{3/2}|K|^2}\right)-\frac{64\kappa}{9}|K|^2\right],
\]
with $\mathcal P$ and $V$ given by
\begin{equation}\label{eq:large-bound}
 \begin{aligned}
 \frac{3(\mathcal P-G)}{2\kappa}
 ={}&\frac{100+12q_0+32q_0^2}{375}
       +\frac{q_0^2}{15+4\kappa}
        \left(\frac{329-256q_0}{1000}+\frac{3}{2n}\right)^2\\
 &+\frac{3q}{8}\left[\frac{q}{3}\left(1+\frac{2}{n}+\frac{4\kappa}{3}\right)
       +\frac{2}{n}-\frac{2}{3}\right],\\
 V(r)={}&\left(1-\frac{2\kappa}{3}\right)\left(\frac{64\kappa}{9}-\frac{8}{3}\right)
       -\frac{4}{n}\left(1+\frac{32\kappa}{9}\right)+\frac{q^2}{6}-\frac{8qr f_3}{3s^{3/2}}\\
 &+\left[\frac{1}{\sqrt{6\kappa}}\left(\frac{37}{12}-\frac{13\kappa}{18}-\frac{6}{n}\right)
       -\frac{8\sqrt\kappa}{3}|r-\beta|\right]_+^2,
 \end{aligned}
\end{equation}
where $\beta=-\frac{3\,\partial_{f_3}F}{s^{3/2}}=(\frac{137}{100}-\frac{\kappa}{15})\frac{f_3}{s^{3/2}}$. For fixed $\kappa,r$, both $q$ and the expression inside the positive part are affine in $\frac{1}{n}$, so $V(r)$ is convex in $\frac{1}{n}$. One half of the second derivative of $\mathcal P$ with respect to $\frac{1}{n}$ is
\[
 \left(\frac{3}{4}-\frac{4\kappa q_0}{3}\right)^2
 +\frac{3\kappa q_0^2}{2(15+4\kappa)}+\frac{8\kappa q_0^2}{n}\ge0 .
\]
It suffices to check $n=5$ and the limit as $n\to\infty$. For $K\ne0$, $\frac{\tr(hK^2)}{\sqrt{s}\,|K|^2}$ is a weighted mean of the $\mu_i$. For $\mu_i<0$, $\frac{f_3}{s^{3/2}}\le(1-\mu_i^2)^{\frac{3}{2}}$, hence $\frac{f_3^2\mu_i^2}{s^3}\le\mu_i^2(1-\mu_i^2)^3\le\frac{27}{256}<\frac{1}{9}$. It therefore follows from \eqref{eq:spectral-parameters} that
\[
 \frac{\tr(h^3K)}{s^{3/2}|K|^2}\ge r_0=\frac{3f_3^2-s^3}{3f_3 s^{3/2}}>0 .
\]

For $\frac{1}{4}\le q_0\le1$, we have $\frac{100+12q_0+32q_0^2}{375}\le0.384$, and completing the square gives $0<\frac{q_0(329-256q_0)}{1000}\le\frac{329^2}{4\cdot256\cdot1000}<0.11$. For $n=5$ and $r\ge r_0$, the squared positive part in $V(r)$ is at most $\frac{5}{18}(\frac{29}{20})^2<\frac{3}{5}$. With $q\le\frac{1}{5}$ and $\frac{64\kappa}{9}>4$, this gives $V(r)\le\frac{3}{2}-\frac{12}{5}+\frac{1}{150}+\frac{3}{5}<0$. Since $\frac{3q_0}{10}\le0.3$ and $\frac{q}{3}(1+\frac{2}{5}+\frac{4\kappa}{3})\le\frac{17}{75}<\frac{4}{15}$, it follows that
\[
 \mathcal P\le\frac{1}{1600}+\frac{2\kappa}{3}
 \left[-0.4+0.384+\frac{0.41^2}{17.4}\right]
 <\frac{1}{1600}-0.004\,\frac{2\kappa}{3}<-10^{-4}.
\]
\begingroup
\rightskip=0pt plus 1fil\relax
As $n\to\infty$, $q=q_0\in[\frac{1}{4},1]$ and $\frac{2\kappa}{3}=\frac{6-4q}{5}$.
We first estimate $\mathcal P+\frac{(52q-1)^2}{25600}$ and then prove $V(r)\le\frac{52q-1}{30}$ for $r\ge r_0$.

 The square term in \eqref{eq:large-bound} satisfies
$\frac{q^2}{15+4\kappa}\left(\frac{329-256q}{1000}\right)^2\le\frac{0.11^2}{17.4}<0.001$.
Adding $\frac{(52q-1)^2}{25600}$ and collecting powers of $1-q$,
we use $(1-q)^j\le(\frac{3}{4})^{j-2}(1-q)^2$ for $j=3,4$ to obtain
\[
 \mathcal P+\frac{(52q-1)^2}{25600}
 \le-0.0043+\frac{23}{5625}(1-q)^2\le-0.002.
\]
\par\endgroup
The maximum of $V$ on $[r_0,\infty)$ occurs at $r_0$ or $\beta$, by convexity up to $\beta$ and monotonicity thereafter. Taking the chord in $q$ gives
\[
 \frac{\beta f_3}{s^{3/2}}=
 \frac{2\kappa}{3}\left(\frac{137}{100}-\frac{\kappa}{15}\right)
 \ge0.53+\frac{74}{75}(1-q)\ge\frac{2\kappa}{3}+0.13.
\]
The formula for $r_0$ gives
\[
 \frac{64\kappa}{9}(\beta-r_0)^2
 \ge\frac{32}{3}\left(0.13+\frac{1}{3}\right)^2
 >2>\frac{1}{6\kappa}\left(\frac{37}{12}-\frac{13\kappa}{18}\right)^2 .
\]
Thus the positive part in $V(r_0)$ vanishes. Evaluating at the two endpoints gives
\[
 \begin{aligned}
 \frac{52q-1}{30}-V(r_0)&\ge\frac{3}{4}-\frac{10}{9}(1-q)(4q-1)\ge\frac{1}{8},\\
 6\kappa\left[\frac{52q-1}{30}-V(\beta)\right]
 &\ge\frac{1}{10}+4(1-q)
 \left[\frac{7}{16}+\left(\frac{8(1-q)}{3}-\frac{3}{4}\right)^2\right]>0 .
 \end{aligned}
\]
Hence $V(r)\le\frac{52q-1}{30}$ for every $r\ge r_0$. At this endpoint, completing the square in $|K|^2$ gives
\[
 \mathcal P+\kappa|K|^2\left(\frac{52q-1}{30}-\frac{64\kappa}{9}|K|^2\right)
 \le\mathcal P+\frac{(52q-1)^2}{25600}\le-0.002 .
\]
\end{proof}

\begin{proof}[Proof of Theorem~\ref{main}]
Suppose $n<S\le\frac{8n}{5}$ and fix the limiting normalized components $\mu_i,T_{ijk}$.
Evaluate $\lambda_i(s)$ and $h_{ijk}(s)$ at $s=\frac{8n}{5}$ and choose the weights in Lemma~\ref{lem:endpoint-verification} for these rescaled principal curvatures.
For these weights, \eqref{eq:exact-endpoint-sign}, \eqref{eq:endpoint-upper-bound} and the lemma give the contradiction
\[
\begin{aligned}
 0\le\Phi_Z\!\left(\frac{8n}{5}\right)
 &\le \mathcal E-\frac{3s}{8|Z|^2}
 \left(\frac{\mathcal C}{|\nabla h|^2}-\frac{f_3}{ns}
       +\frac{32\kappa\,\tr(h^3K)}{9s}\right)^2\\
 &\quad-\frac{64\kappa^2}{9}
 \left(\frac{|K|^2}{2}-\frac{|hK|^2}{s}\right)<-10^{-4}.
\end{aligned}
\]
\end{proof}

\begin{rem}
In fact, we can improve the bound $S>\frac{8n}{5}$ to the stronger bound $S>\frac{32}{19}n$ via a more involved argument.
\end{rem}

\end{document}